\documentclass[12pt,twoside]{amsart}
\usepackage{latexsym,amssymb,amsmath}
\numberwithin{equation}{section}
\newtheorem{theorem}{Theorem}[section]
\newtheorem{lemma}[theorem]{Lemma}

\newtheorem{corollary}[theorem]{Corollary}

\theoremstyle{definition}
\newtheorem{definition}[theorem]{Definition} % \theoremstyle{remark}
\newtheorem{remark}[theorem]{Remark}
\newtheorem{example}[theorem]{Example}

\begin{document}

%%%%%%%%%%%%%%%%%%%%%%%%%%%%%%%%%%%%%%%%%%%%%%%%%%%%%%%%%%%%%%%%%%
% New Commands

\newcommand{\m}[1]{\marginpar{\addtolength{\baselineskip}{-3pt}{\footnotesize
\it #1}}} 
\newcommand{\A}{\mathcal{A}}
\newcommand{\K}{\mathcal{K}} 
\newcommand{\knd}{\mathcal{K}^{[d]}_n}
\newcommand{\F}{\mathcal{F}}
\newcommand{\N}{\mathbb{N}}
\newcommand{\pr}{\mathbb{P}}
\newcommand{\I}{\mathit{I}}
\newcommand{\G}{\mathcal{G}}
\newcommand{\D}{\mathcal{D}}
\newcommand{\lcm}{\operatorname{lcm}}
\newcommand{\ndp}{N_{d,p}}
\newcommand{\tor}{\operatorname{Tor}}
\newcommand{\reg}{\operatorname{reg}} 
\newcommand{\mf}{\mathfrak{m}}

\def\bb{{{\rm \bf b}}}
\def\cc{{{\rm \bf c}}}

%%%%%%%%%%%%%%%%%%%%%%%%%%%%%%%%%%%%%%%%%%%%%%%%%%%%%%%%%%%%%%%% 
 
\title{Shellable graphs and 
sequentially Cohen-Macaulay bipartite graphs}
\thanks{Version: \today}
 
\author{Adam Van Tuyl}
\address{Department of Mathematical Sciences \\
Lakehead University \\
Thunder Bay, ON P7B 5E1, Canada}
\email{avantuyl@sleet.lakeheadu.ca}
\urladdr{http://flash.lakeheadu.ca/$\sim$avantuyl/}

\author{Rafael H. Villarreal}
\address{
Departamento de
Matem\'aticas\\
Centro de Investigaci\'on y de Estudios
Avanzados del
IPN\\
Apartado Postal
14--740 \\
07000 Mexico City, D.F.
}
\email{vila@math.cinvestav.mx}
\urladdr{http://www.math.cinvestav.mx/$\sim$vila/}
 

\begin{abstract}
Associated to a simple undirected graph $G$ is a simplicial complex 
$\Delta_G$ whose faces correspond
to the independent sets of $G$.  We call a graph $G$ 
shellable if $\Delta_G$ is a shellable simplicial complex in the non-pure
sense of Bj\"orner-Wachs.  We are then interested in
determining what families of graphs have the property that $G$ is
shellable.  We show that all chordal graphs are shellable.  Furthermore, 
we classify all the shellable bipartite graphs;  they are precisely
the sequentially Cohen-Macaulay 
bipartite graphs. 
We also give a recursive procedure to verify if a bipartite
graph is shellable.
Because shellable implies that the associated Stanley-Reisner ring is
sequentially Cohen-Macaulay, our results complement and extend recent 
work on the problem of determining when the edge ideal of a graph is
(sequentially) Cohen-Macaulay.  We also give a new proof
for a result of Faridi on the sequentially Cohen-Macaulayness of
simplicial forests. 
\end{abstract}
 
\maketitle

%%%%%%%%%%%%%%%%%%%%%%%%%%%%%%%%%%%%%%%%%%%%%%%%%%%%%%%%%%%%%%%%%%%

\section{Introduction} 

Let $G$ be a simple (no loops or multiple edges) undirected graph
on the vertex set $V_G = \{x_1,\ldots,x_n\}$.   By identifying the vertex $x_i$
with the variable $x_i$ in the polynomial ring $R = k[x_1,\ldots,x_n]$
over a field $k$, we can
associate to $G$ a quadratic square-free monomial ideal
$I(G) = ( \{x_ix_j ~|~ \{x_i,x_j\} \in E_G\})$ where $E_G$ is the edge set of $G$.
The ideal $I(G)$ is called the {\it edge ideal} of $G$.
Using the Stanley-Reisner correspondence,
we can associate to $G$ the simplicial complex $\Delta_G$ where
$I_{\Delta_G} = I(G)$.   
Notice that the faces of $\Delta_G$ are the {\it independent
sets} or {\it stable sets\/} of $G$. Thus $F$ is a face of $\Delta_G$ 
if and only if there is no edge of $G$ joining any two vertices of
$F$. The dual concept of an independent set is a
{\it vertex cover\/}, i.e., a subset $C$ of $V_G$ is a vertex cover of
$G$  if and only if $V_G\setminus C$ is an independent set of $G$. 

We call a graph $G$ {\it {\rm (}sequentially{\rm )} Cohen-Macaulay}
if $R/I(G)$ is (sequentially) Cohen-Macaulay. 
Recently, a number of authors (for example, see
\cite{EV,FH,FVT2,herzog-hibi,HHZ,Vi2,unmixed}) have been interested in  
classifying or identifying (sequentially) Cohen-Macaulay graphs $G$ in terms of 
the combinatorial properties of $G$.   This paper complements and extends some
of this work by introducing the notion of a {\it shellable graph}.
We shall call a graph $G$ shellable if the simplicial complex $\Delta_G$
is a shellable simplicial complex 
(see Definition \ref{shellabledefn}).  Here, we mean
the non-pure definition of shellability as introduced by Bj\"orner
and Wachs \cite{BW}.   
Because a shellable simplicial complex has the property that its associated 
Stanley-Reisner ring is sequentially Cohen-Macaulay, by identifying
shellable graphs, we are in fact identifying some
of the sequentially Cohen-Macaulay graphs.

We begin in Section 2 by formally introducing shellable graphs and
discussing some of their basic properties.  We then focus on the shellability
of bipartite graphs.
 Recall that a graph $G$ is {\it bipartite\/} if the vertex set $V_G$
can be partitioned into two disjoint sets $V  = V_1 \cup V_2$ such
that every edge of $G$ contains one vertex in $V_1$ and the other in $V_2$.
Furthermore, let $N_G(x)$ denote the set of {\it neighbors} of the
vertex $x$.
We then show:

\begin{theorem}[Corollary \ref{recursivebuild}]\label{thm1}
Let $G$ be a bipartite graph.  Then $G$ is shellable if and only
if there are adjacent vertices $x$ and $y$ 
with $\deg(x)=1$ such that 
the bipartite graphs $G \setminus (\{x\} \cup N_G(x))$
and $G \setminus (\{y\} \cup N_G(y))$ are shellable.
\end{theorem}

We also consider the shellability of chordal graphs.
A graph ${G}$ is  
{\it chordal\/} (or {\it triangulated}) if every cycle 
${\mathcal C}_n$ of ${G}$ of length $n\geq 4$ has a 
chord. A {\it chord} of ${\mathcal C}_n$ is an edge joining two 
non-adjacent vertices of ${\mathcal C}_n$.   Chordal graphs
then have a nice combinatorial property:

\begin{theorem}[Theorem \ref{chordaltheorem}]
Let $G$ be a chordal graph.  Then $G$ is shellable.
\end{theorem}

\noindent
Because $G$ being shellable implies that $G$ is sequentially 
Cohen-Macaulay, the above result gives a new proof to the main result
of Francisco and the first author \cite{FVT2} that all chordal
graphs are sequentially Cohen-Macaulay.  

The main result of Section 3 is to classify all sequentially Cohen-Macaulay
bipartite graphs.  Precisely, we show:

\begin{theorem}[Theorem \ref{SCM=shellable}] \label{them3}
Let $G$ be a bipartite graph.
Then $G$ is sequentially Cohen-Macaulay if and only $G$ is
shellable.
\end{theorem}

\noindent
Note that all shellable graphs are automatically sequentially Cohen-Macaulay
(see Stanley \cite{Stanley} or Theorem \ref{shellable->scm}), but
the converse is not true in general.  So,
the above theorem says that among the bipartite graphs, those that are sequentially
Cohen-Macaulay are precisely those that are shellable.  This generalizes
a result of Estrada and the second author \cite{EV} which
showed that $G$ is a Cohen-Macaulay bipartite graph if and only
if $\Delta_G$ has a pure shelling.  
Because we can use Theorem \ref{thm1} to recursively check if a 
bipartite graph is shellable, Theorem \ref{them3} implies we can
verify recursively if a 
bipartite graph is sequentially Cohen-Macaulay. 

In the fourth section we consider connected bipartite graphs with 
bipartition $V_1=\{x_1,\ldots,x_g\}$ and $V_2=\{y_1,\ldots,y_g\}$ such
that $\{x_i,y_i\}\in E_G$ for all $i$ and $g\geq 2$. 
Following Carr\'a Ferro and Ferrarello \cite{carra-ferrarello},
we can associate to $G$ a directed graph $\mathcal{D}$.  Carr\'a
Ferro and Ferrarello 
gave an alternative classification of Cohen-Macaulay bipartite graphs in terms
of the properties of $\mathcal{D}$ (the original classification is 
due of Herzog and Hibi \cite{herzog-hibi}).  We show
how $G$ being sequentially Cohen-Macaulay affects the graph $\mathcal{D}$. 

In the final section we extend the scope of our investigation
to include the edge ideals associated to clutters (a type of hypergraph).  
As in the graph case, we say that a clutter
$\mathcal{C}$ is shellable if the simplicial complex associated to
the edge ideal $I(\mathcal{C})$ 
is a shellable simplicial complex.    
We show  (the free vertex property is defined in Section 5):
\begin{theorem}[Theorem \ref{fvp}]
If a clutter $\mathcal{C}$ has the free vertex property, then
$\mathcal{C}$ is shellable.  
\end{theorem} 
By applying
a result of Herzog, Hibi, Trung and Zheng \cite{hhtz}, we recover as a corollary
the fact that all simplicial forests are sequentially Cohen-Macaulay.
This result was first proved by Faridi \cite{Faridi}.

%%%%%%%%%%%%%%%%%%%%%%%%%%%%%%%%%%%%%%%%%%%%%%%%%%%%%%%%%%%%%%%%%%%

\section{Shellable graphs}

We continue to use the notation and definitions used in the introduction.  In
this section we introduce shellable graphs, describe some
of their properties, and identify some families of shellable graphs.

\begin{definition}\label{shellabledefn}
A simplicial complex $\Delta$ is {\it shellable\/}  if 
the facets (maximal faces) of $\Delta$ can be ordered $F_1,\ldots,F_s$ such that 
for all $1\leq i<j\leq s$, there 
exists some $v\in F_j\setminus F_i$ and some 
$\ell\in \{1,\ldots,j-1\}$ with $F_j\setminus F_\ell= \{v\}$. We call
$F_1,\ldots,F_s$ a {\it shelling} of $\Delta$ when the facets have been
ordered with respect to the shellable definition.  For a fixed
shelling of $\Delta$, if $F,F' \in \Delta$
then we write $F < F'$ to mean that $F$
appears before $F'$ in the ordering.
\end{definition}

\begin{remark}
The above definition of shellable is due to
Bj\"orner and Wachs \cite{BW} and is usually referred
to as {\it nonpure shellable}, although
in this paper we will drop the adjective ``nonpure''.  Originally, the definition of
shellable also required that 
the simplicial complex be pure, that is, all the facets have same dimension.
We will say $\Delta$ is {\it pure shellable} if it also satisfies this
hypothesis.
\end{remark}

\begin{definition} Let $G$ be a simple undirected graph with associated
simplicial complex $\Delta_G$.  We say
$G$ is a {\it shellable graph} if $\Delta_G$ is a shellable simplicial 
complex.
\end{definition}

To prove that a graph $G$ is shellable, it suffices to prove each
connected component of $G$ is shellable, as demonstrated below.

\begin{lemma}\label{oct14-06} Let $G_1$ and $G_2$ be two graphs with
disjoint sets of 
vertices and let $G=G_1\cup G_2$. 
Then $G_1$ and $G_2$ are shellable if and only if $G$ 
is shellable.  
\end{lemma}

\begin{proof} $(\Rightarrow)$
 Let $F_{1},\ldots,F_{r}$ and  $H_{1},\ldots,H_{s}$ be 
the shellings of $\Delta_{G_1}$ and $\Delta_{G_2}$ respectively. 
Then
if we order the facets of $\Delta_G$ as
$$
F_{1}\cup H_{1},\ldots,F_{1}\cup H_{s};\, 
F_{2}\cup H_{1},\ldots,F_{2}\cup H_{s};\, \ldots;\, 
F_{r}\cup H_{1},\ldots,F_{r}\cup H_{s}
$$
we get a shelling of $\Delta_G$. Indeed if $F'<F$ are two facets
of $\Delta_G$ we have two cases to consider. 
Case (i): $F'=F_i\cup H_k$ and $F=F_j\cup H_t$, where $i<j$.  Because
$\Delta_{G_1}$ is shellable there is $v\in F_j\setminus F_i$ and 
$\ell<j$ with $F_j\setminus F_\ell= \{v\}$. 
Hence $v\in F\setminus F'$, $F_\ell\cup H_t<F$, and 
$F\setminus(F_\ell\cup H_t)=\{v\}$. Case (ii): $F'=F_k\cup H_i$ and
$F=F_k\cup H_j$, where $i<j$. This case follows from the 
shellability of $\Delta_{G_2}$.  

$(\Leftarrow)$  Note that if $F$ is a facet of $\Delta_G$,
then $F' = F \cap V_{G_1}$, respectively, $F'' = F \cap V_{G_2}$, is a facet
of $\Delta_{G_1}$, respectively $\Delta_{G_2}$.  We now show that
$G_1$ is shellable and omit the similar proof for the shellability of $G_2$.
Let $F_1,\ldots,F_t$ be a shelling of $\Delta_G$, and consider the subsequence
\[F_{i_1},\ldots,F_{i_s} ~~~\mbox{with $1 = i_1 < i_2 < \cdots < i_s$} \]
where $F_1 \cap V_{G_2} = F_{i_j} \cap V_{G_2}$ for $i_j \in \{i_1,\ldots,i_s\}$,
but $F_1 \cap V_{G_2} \neq F_k \cap V_{G_2}$ for any $k \in \{1,\ldots,t\} \setminus 
\{i_1,\ldots,i_s\}$.  We then claim that 
\[F'_1 = F_{i_1} \setminus V_{G_2},\ F'_2 = F_{i_2} \setminus V_{G_2},
\ldots, F'_s = F_{i_s} \setminus V_{G_2}\]
is a shelling of $\Delta_{G_1}$.  We first
show that this is a complete list of facets; indeed, each $F'_j = F_{i_j} \cap V_{G_1}$
is a facet of $\Delta_{G_1}$, and furthermore, for any facet
$F \in \Delta_{G_1}$, $F \cup (F_1 \cap V_{G_2})$ is a facet of $\Delta_G$,
and hence $F \cup (F_1 \cap V_{G_2}) = F_{i_j}$ for some $i_j \in \{i_1,\ldots,i_s\}$.

Because the $F_i$'s form a shelling, if $1 \leq k < j \leq s$,
there exists $v \in F_{i_j} \setminus F_{i_k} = (F_{i_j} \setminus V_{G_2})\setminus
(F_{i_k}\setminus V_{G_2}) = F'_j \setminus F'_k$ such that
$\{v\} = F_{i_j} \setminus F_{\ell}$ for some $1 \leq \ell < i_j$.
It suffices to show that $F_{\ell}$ is among $F_{i_1},\ldots,F_{i_s}$.  Now
because $F_{i_j} \cap V_{G_2} \subset F_{i_j}$ and $v \not\in
F_{i_j} \cap V_{G_2}$, we must have $F_{i_j} \cap V_{G_2} \subset F_{\ell}$.
So, $F_{\ell} \cap V_{G_2} \supset F_{i_j} \cap V_{G_2}$. 
But $F_{\ell} \cap V_{G_2}$ is a facet of $\Delta_{G_2}$,
so we must have $F_{\ell} \cap V_{G_2} = F_{i_j} \cap V_{G_2}$.
So  $F_{\ell} = F_{i_r}$ for some $r < j$,
and hence, $\{v\} = F'_j \setminus F'_r$, as desired. 
\end{proof}

Given a subset $S \subset V_G$, by $G \setminus S$, we mean
the graph formed from $G$ by deleting all the vertices in $S$, and all
edges incident to a vertex in $S$. If $x$ is a vertex of $G$, then
its {\it neighbor set},  
denoted by $N_G(x)$, is
the set of vertices of $G$ adjacent to $x$.   
If $F$ is a face of a simplicial complex $\Delta$, the {\it link} of $F$ is
defined to be $\operatorname{lk}_{\Delta}(F) = \{G ~|~ G \cup F \in
\Delta,~~ G \cap F = \emptyset\}$. 
When $F = \{x\}$, then we shall abuse notation and
write $\operatorname{lk}_{\Delta}(x)$ 
instead of $\operatorname{lk}_{\Delta}(\{x\})$.

\begin{lemma}\label{link}
 Let $x$ be a vertex of $G$ and let
$G' = G\setminus(\{x\}\cup N_G(x))$.  Then
\[\Delta_{G'} = \operatorname{lk}_{\Delta_G}(x).\]
In particular, $F$ is a facet of $\Delta_{G'}$ if and only if 
$x\notin F$ and $F \cup
\{x\}$ is a facet 
of $\Delta_G$.
\end{lemma}

\begin{proof}
If $F \in \operatorname{lk}_{\Delta_G}(x)$,
then $x \not\in F$, and $F \cup \{x\} \in \Delta_G$ implies that $F \cup \{x\}$
is an independent set of $G$.  So $(F \cup \{x\}) \cap N_G(x) = \emptyset$.
But this means that $F \subset V_{G'}$ because  
$V_{G'}= V_G \setminus (\{x\} \cup N_G(x))$.
Thus $F \in \Delta_{G'}$ since $F$ is also
an independent set of the smaller graph $G'$.

Conversely, if $F \in \Delta_{G'}$, then $F$ is an independent set of $G'$
that does not contain any of the vertices of $\{x\} \cup N_G(x)$.  But
then $F \cup \{x\}$ is an independent set of $G$, i.e., $F \cup \{x\} \in \Delta_G$.
So $F \in \operatorname{lk}_{\Delta_G}(x)$.

The last statement follows readily from the fact that $F$ is a facet
of $\operatorname{lk}_{\Delta_G}(x)$ if and only if $x\notin F$ and 
$F \cup \{x\}$ is a facet of $\Delta_G$. 
\end{proof}

The property of shellability is preserved when
removing the vertices $\{x\} \cup N_G(x)$ and all incident
edges from $G$ for
any vertex $x$.

\begin{theorem} \label{removevertex-shellable}
 Let $x$ be a vertex of $G$ and let
$G' = G\setminus(\{x\}\cup N_G(x))$.
If $G$ is shellable, then $G'$ is shellable.
\end{theorem}

\begin{proof}
Let $F_1,\ldots,F_s$ be a shelling of $\Delta_G$.  Suppose the subsequence 
\[F_{i_1},F_{i_2},\ldots,F_{i_t} ~~\mbox{with $i_1 < i_2 < \cdots < i_t$}\]
is the list of all the facets with $x \in F_{i_j}$.  Setting $H_j =
F_{i_j}\setminus\{x\}$ 
for each $j =1,\ldots,t$, Lemma \ref{link} implies that the $H_j$'s
are the facets of $\Delta_{G'}$. 

We claim that $H_1,\ldots,H_t$ is a shelling of $\Delta_{G'}$.  Because
the $F_i$'s form a shelling, if $1 \leq k < j \leq t$, there
exists a vertex $v \in F_{i_j} \setminus F_{i_k} = (F_{i_j}
\setminus\{x\}) \setminus 
(F_{i_k}\setminus\{x\}) = 
(H_j \setminus H_k)$ such that $\{v\} = F_{i_j} \setminus F_{\ell}$
for some $1 \leq \ell < i_j$.   It suffices to show that $F_\ell$ is
among the list $F_{i_1}, 
\ldots, F_{i_t}$.  But because $x \in F_{i_j}$ and $x \neq v$, we must
have $x \in F_{\ell}$.  Thus $F_{\ell} = F_{i_k}$ for some $k \leq j$.  But
then $\{v\} = F_{i_j} \setminus F_{\ell} = H_{j} \setminus H_k$.  So,
the $H_i$'s form 
a shelling of $\Delta_{G'}$.
\end{proof}

Let $G$ be a graph and let $S\subset V_G$. For use below consider the
graph $G\cup W_G(S)$ obtained from $G$ by
adding new vertices $\{y_i~\vert\,~ x_i\in S\}$ and new
edges $\{\{x_i,y_i\}~\vert\,~ x_i\in S\}$. The edges $\{x_i,y_i\}$ are
called {\it whiskers}.   The notion of a whisker was introduced by
the second author \cite{ITG,Vi2} to study how modifying the graph $G$ 
affected the Cohen-Macaulayness of $G$; this
idea was later generalized by Francisco and H\`a \cite{FH} in
their study of sequentially Cohen-Macaulay graphs.  We can give
a shellable analog of \cite[Theorem 4.1]{FH}.

\begin{corollary} Let $G$ be a graph and let 
$S\subset V_G$. If $G\cup W_G(S)$ is shellable, then 
$G\setminus S$ is shellable. 
\end{corollary}

\begin{proof} We may assume that $S=\{x_1,\ldots,x_s\}$. Set 
$G_0=G \cup W_G(S)$ and $G_i=G_{i-1}\setminus(\{y_i\}\cup N_G(y_i))$ for
$i=1,\ldots,s$. Notice that $G_s=G\setminus S$. Hence, by repeatedly applying
Theorem \ref{removevertex-shellable}, the graph $G\setminus S$ is 
shellable. \end{proof}

We now turn our attention to the shellability of bipartite graphs.

\begin{lemma} \label{degree1} Let $G$ be a bipartite graph with 
bipartition $\{x_1,\ldots,x_m\}$, $\{y_1,\ldots,y_n\}$. If $G$
is shellable and $G$ has no isolated vertices, 
then there is $v\in V_G$ with $\deg(v)=1$.
\end{lemma}

\begin{proof} Let $F_1,\ldots,F_s$ be a shelling of $\Delta_G$. We may
assume that $F_i=\{y_1,\ldots,y_n\}$, $F_j=\{x_1,\ldots,x_m\}$ and
$i<j$. Then there is $x_k\in F_j\setminus F_i$ and $F_\ell$ with
$\ell\leq j-1$ such that $F_j\setminus F_\ell=\{x_k\}$. For simplicity
assume that $x_k=x_1$. Then $\{x_2,\ldots,x_m\}\subset F_\ell$ and
there is $y_t$ in $F_\ell$ for some $1\leq t\leq n$. Since
$\{y_t,x_2,\ldots,x_m\}$ is an independent set of $G$, we get that
$y_t$ can only be adjacent to $x_1$. Thus $\deg(y_t)=1$ because
$G$ has no isolated vertices. \end{proof}

\begin{theorem}\label{oct15-1-06} Let $G$ be a graph and let $x_1,y_1$ be two
adjacent vertices of $G$ with $\deg(x_1)=1$. If
$G_1=G\setminus(\{x_1\}\cup N_G(x_1))$ and
$G_2=G\setminus(\{y_1\}\cup N_G(y_1))$, then $G$ is shellable
if and only if $G_1$ and $G_2$ are shellable. 
\end{theorem}

\begin{proof} If $G$ is shellable, then $G_1$ and $G_2$ are shellable by
Theorem \ref{removevertex-shellable}.  So it suffices to prove the
reverse direction.
Let $F_1',\ldots,F_r'$ be a shelling of 
$\Delta_{G_1}$ and let $H_1',\ldots,H_s'$ be a shelling of 
$\Delta_{G_2}$. It suffices to prove that
$$
F_1'\cup\{x_1\},\ldots,F_r'\cup\{x_1\},H_1'\cup\{y_1\},\ldots,
H_s'\cup\{y_1\}
$$
is a shelling of $\Delta_G$.  One first shows that this is the complete list
of facets of $\Delta_G$ using Lemma~\ref{link}.  Indeed,
take any facet $F$ of $\Delta_G$. If $y_1 \in F$,
then $x_1 \not\in F$ because $\{x_1,y_1\}$ is an edge of $G$,
and by Lemma ~\ref{link}, $F \setminus \{y_1\} = H_i'$ for
some $i$.  On the other hand, if $y_1 \not\in F$,
we must have $x_1 \in F$, because if not, then
$\{x_1\} \cup F$ is larger independent set of $G$ because
$x_1$ is only adjacent to $y_1$.  Again, by Lemma ~\ref{link},
we have $F \setminus \{x\} = F_i'$ for some $i$.  
Let $F'<F$ be two facets of $\Delta_G$. 
There are three cases to consider. Case (i): $F'=F_i'\cup\{x_1\}$ 
and $F=H_j'\cup\{y_1\}$. Since $H_j'\cup\{x_1\}$ is an independent
set of $G$, it is contained in a facet of $\Delta_G$, i.e., 
$H_j'\cup\{x_1\}\subset F_\ell'\cup\{x_1\}$ for some $\ell$. Hence
$(H_j'\cup\{y_1\})\setminus(F_\ell'\cup\{x_1\})=\{y_1\}$,
$y_1\in F\setminus F'$, and $F_\ell'\cup\{x_1\}<F$. The remaining two
cases follow readily from the shellability of $\Delta_{G_1}$ 
and $\Delta_{G_2}$. \end{proof}

Putting together the last two results yields 
a recursive procedure to verify if a bipartite
graph is shellable.

\begin{corollary} \label{recursivebuild}
Let $G$ be a bipartite graph.  Then $G$ is shellable if and only
if there are adjacent vertices $x$ and $y$ 
with $\deg(x)=1$ such that 
the bipartite graphs $G \setminus (\{x\} \cup N_G(x))$
and $G \setminus (\{y\} \cup N_G(y))$ are shellable.
\end{corollary}

\begin{proof}  By Lemma \ref{oct14-06} it suffices to
verify the statement when $G$ is connected.  
By Lemma \ref{degree1} there exists
a vertex of $x_1$ with $\deg(x_1) =1$.   Now apply the previous
theorem.
\end{proof}

\begin{example}  
The {\it complete bipartite graph}, denoted $\mathcal{K}_{m,n}$,
is the graph with vertex set $V_G = \{x_1,\ldots,x_m,y_1,\ldots,y_n\}$
and edge set $E_G = \{\{x_i,y_j\} ~|~ 1 \leq i \leq m, 1 \leq j \leq n\}$.
If $m,n \geq 2$, then $\mathcal{K}_{m,n}$
is not shellable since the graph has no vertex of degree one.
On the other hand, if $m =1$ and $n \geq 1$, then
$\mathcal{K}_{m,n}$ is shellable since 
the only facets are $F_1 = \{y_1,\ldots,y_n\}$
and $F_2 = \{x_1\}$ and we have a shelling with $F_1 < F_2$.  Similarly,
$\mathcal{K}_{m,1}$ is shellable for all $m \geq 1$.  
\end{example}

We complete this section by showing that all chordal
graphs are shellable.
A graph ${G}$ is {\it triangulated\/} or 
{\it chordal\/} if every cycle ${\mathcal C}_n$ of ${G}$ of length $n\geq 4$ has a 
chord. A {\it chord} of ${\mathcal C}_n$ is an edge joining two 
non-adjacent vertices of ${\mathcal C}_n$. Let $S$ be a set of
vertices of a graph $G$.  The {\it induced
subgraph\/} $G_S$ is the maximal subgraph of $G$ with
vertex set $S$.  
For use below we call a complete subgraph of $G$ a {\it clique}. As usual, 
a complete graph with $r$ vertices is denoted by ${\mathcal K}_r$.

\begin{lemma}{\rm\cite[Theorem~8.3]{Toft}}\label{toft-lemma}
Let $G$ be a chordal graph and let ${\mathcal K}$ be a complete 
subgraph of $G$. If ${\mathcal K}\neq G$, then there is $x\not\in
V({\mathcal K})$ 
such that $G_{N_G(x)}$ is a complete subgraph.
\end{lemma}

\begin{theorem} \label{chordaltheorem}
Let $G$ be a chordal graph. Then $G$ is 
shellable.
\end{theorem}

\begin{proof} We proceed by induction on $n = |V_G|$. Let 
$V_G=\{x_1,\ldots,x_n\}$ be the vertex set of $G$. If $G$ is a
complete graph, then $\Delta_G$ consists of $n$ isolated vertices and
they clearly form a shelling. Thus by Lemma~\ref{oct14-06} we may assume
that $G$ is a connected non-complete graph. According to
Lemma~\ref{toft-lemma} there is  
$x_1\in V_G$ such that $G_{N_G(x_1)}={\mathcal K}_{r-1}$ is a
complete subgraph for some $r \geq 1$.  (To apply Lemma ~\ref{toft-lemma}, take
$\mathcal K$ to be any edge of $G$;  this is clearly a complete graph.)
Notice that 
$G_{\{x_1\}\cup N_G(x_1)} ={\mathcal K}_r$ and that ${\mathcal
K}_r$ is the 
only maximal complete subgraph of $G$ that contains $x_1$. We may
assume that $V({\mathcal K}_r)=\{x_1,\ldots,x_r\}$.
Consider the subgraphs $G_i=G\setminus(\{x_i\}\cup N_G(x_i))$,
which are also chordal.   By
induction there is a shelling $F_{i1},\ldots,F_{is_i}$ 
of $\Delta_{G_i}$ for $i=1,\ldots,r$. Observe that any facet of
$\Delta_G$ intersects $V({\mathcal K}_r)$ in exactly one vertex. Thus
by Lemma~\ref{link} the following is the 
complete list of facets of $\Delta_G$:
$$
F_{11}\cup\{x_1\},\ldots,F_{1s_1}\cup\{x_1\};
F_{21}\cup\{x_2\},\ldots,F_{2s_2}\cup\{x_2\};\ldots;
F_{r1}\cup\{x_r\},\ldots,F_{rs_r}\cup\{x_r\}.
$$
We claim that this linear ordering is a shelling of $\Delta_G$. Let
$F'<F$ be two facets of $\Delta_G$. There are two cases to consider.
Case (i): $F'=F_{ik}\cup\{x_i\}$ and $F=F_{jt}\cup\{x_j\}$, where
$i<j$. Notice that $F_{jt}\cup\{x_1\}$ is an independent set of $G$ 
because $F_{jt}\cap V(\mathcal{K}_r)=\emptyset$. Thus
$F_{jt}\cup\{x_1\}$ can be extended to a facet of $G$, i.e., 
$F_{jt}\cup\{x_1\}\subset F_{1\ell}\cup\{x_1\}$ for some
$1\leq\ell\leq s_1$. Set $F''= F_{1\ell}\cup\{x_1\}$. Hence $x_j\in
F\setminus F'$, $F\setminus F''=\{x_j\}$, and $F''<F$. Case (ii): 
$F'=F_{ik}\cup\{x_i\}$ and $F=F_{it}\cup\{x_i\}$, with $k<t$. 
This case follows from the 
shellability of $\Delta_{G_i}$. \end{proof}

\begin{remark}
As shown below (Theorem \ref{shellable->scm}), if a graph $G$
is shellable, then it is also sequentially Cohen-Macaulay.  The above
theorem, therefore, gives a new proof to the fact that all chordal
graphs are sequentially Cohen-Macaulay as first proved in \cite{FVT2}.
To show that all chordal graphs are 
sequentially Cohen-Macaulay,  the authors of \cite{FVT2} 
show that for each degree $d \geq 0$,  the square-free part
of the Alexander dual $I(G)^{\vee}$ (also
defined below) in degree $d$ has {\it linear quotients},
that is, there is an ordering of the generators $\{u_1,\ldots,u_s\}$
of the square-free part of $I(G)^{\vee}$ of degree $d$
such that $(u_1,\ldots,u_{i-1}):(u_i) = (x_{i_1},\ldots,x_{i_t})$
for $i=1,\ldots,s$.
However, when $G$ is shellable, the 
generators of the Alexander dual $I(G)^{\vee}$ 
must have linear quotients (see \cite[Theorem 1.4(c)]{hhz-ejc} and 
\cite{SJZ});  
so, when $G$ is chordal,
the ideal $I(G)^{\vee}$ also has linear quotients, a fact, to the
best of our knowledge, that has never been noticed.
\end{remark}

%%%%%%%%%%%%%%%%%%%%%%%%%%%%%%%%%%%%%%%%%%%%%%%%%%%%%%%%%%%%%%%%%%%

\section{Sequentially Cohen-Macaulay bipartite graphs}

In this section we classify all sequentially Cohen-Macaulay bipartite graphs.
We begin by recalling the relevant definitions and results about sequentially
Cohen-Macaulay modules.

\begin{definition} \label{d.seqcm}
Let $R=k[x_1,\dots,x_n]$. A graded $R$-module $M$ is called 
{\it sequentially Cohen-Macaulay} (over $k$)
if there exists a finite filtration of graded $R$-modules
\[ 0 = M_0 \subset M_1 \subset \cdots \subset M_r = M \]
such that each $M_i/M_{i-1}$ is Cohen-Macaulay, and the Krull dimensions of the
quotients are increasing:
\[\dim (M_1/M_0) < \dim (M_2/M_1) < \cdots < \dim (M_r/M_{r-1}).\]
\end{definition}
 
As first shown by Stanley \cite{Stanley}, shellable implies
sequentially Cohen-Macaulay. 

\begin{theorem}\label{shellable->scm}
Let $\Delta$ be a simplicial complex, and suppose that $R/I_{\Delta}$
is the associated Stanley-Reisner ring.  If $\Delta$ is shellable, then
$R/I_{\Delta}$ is sequentially Cohen-Macaulay.  
\end{theorem}

We now specialize to the case of graphs by providing a sequentially
Cohen-Macaulay analog of Theorem \ref{removevertex-shellable}.

\begin{theorem} \label{gscm->cscm1}
 Let $x$ be a vertex of $G$ and let
$G' = G\setminus(\{x\}\cup N_G(x))$.
If $G$ is sequentially Cohen-Macaulay, then $G'$ is sequentially Cohen-Macaulay.

\end{theorem}

\begin{proof}
Let $F_1,\ldots,F_s$ be the facets of
$\Delta=\Delta_G$, and let $F_1,\ldots,F_r$ be the facets of $\Delta$ that
contain $x$.   Set $\Gamma=\Delta_{G'}$;  by Lemma~\ref{link}, the
facets of $\Gamma$ are
$F_1'=F_1\setminus\{x\},\ldots,F_r'=F_r\setminus\{x\}$.

Consider the pure simplicial complexes 
\begin{eqnarray*}
\Delta^{[k]}&=&\langle\{F\in\Delta\vert\, \dim(F)=k\}\rangle;\ \ -1\leq
k\leq\dim(\Delta),\\
\Gamma^{[k]}&=&\langle\{F\in\Gamma\vert\, \dim(F)=k\}\rangle;\ \ -1\leq
k\leq\dim(\Gamma),
\end{eqnarray*}
where $\langle{\mathcal F}\rangle$ denotes the subcomplex generated by 
the set of faces $\mathcal F$. Recall that $H$ is a face of 
$\langle{\mathcal F}\rangle$ if and only if $H$ is contained in 
some  $F$ in $\mathcal{F}$. Take a facet $F_i'$ of $\Gamma$ of
dimension $d=\dim(\Gamma)$. Then
$F_i'\cup\{x\}\in\Delta^{[d+1]}$ and consequently $\{x\}\in
\Delta^{[k+1]}$ for $k\leq d$. Because the 
facets of $\Gamma$ are
$F_1'=F_1\setminus\{x\},\ldots,F_r'=F_r\setminus\{x\}$,  we have
the equality  
$$
\Gamma^{[k]}={\rm lk}_{\Delta^{[k+1]}}(x) 
$$
for $k\leq d$. By \cite[Theorem~3.3]{duval} the complex
$\Delta$ is sequentially Cohen-Macaulay if and only if $\Delta^{[k]}$
is Cohen-Macaulay for  
$-1\leq k\leq \dim(\Delta)$.  Because $\Delta^{[k]}$ is Cohen-Macaulay,
by \cite[Proposition 5.3.8]{V} 
${\rm lk}_{\Delta^{[k]}}(F)$ is Cohen-Macaulay for any $F \in
\Delta^{[k]}$. Thus, $\Gamma^{[k]}={\rm lk}_{\Delta^{[k+1]}}(x)$ is Cohen-Macaulay
 for any $-1 \leq k \leq \dim (\Gamma) \leq \dim(\Delta)-1$.
Therefore $\Gamma$ is sequentially Cohen-Macaulay by
\cite[Theorem~3.3]{duval}, as
required. \end{proof} 

\begin{example}  The six cycle $\mathcal{C}_6$ is a counterexample to the
converse of the above statement.  For any vertex $x$ of $\mathcal{C}_6$,
the graph $\mathcal{C}_6 \setminus (\{x\}\cup N_G(x))$ is a tree,
which is sequentially Cohen-Macaulay.  (A tree is a chordal
graph, so by Theorem \ref{chordaltheorem}, a tree is shellable, and hence,
sequentially Cohen-Macaulay
by Theorem \ref{shellable->scm}.)  However, the only
sequentially Cohen-Macaulay cycles are $\mathcal{C}_3$ and $\mathcal{C}_5$
\cite[Proposition 4.1]{FVT2}.
\end{example}

A corollary of the above result is the following result of 
Francisco and H\`a.  Here $W_G(S)$
is the whisker notation introduced in the previous section.

\begin{corollary}{\rm \cite[Theorem~4.1]{FH}} Let $G$ be a graph and let 
$S\subset V_G$. If $G\cup W_G(S)$ is sequentially Cohen-Macaulay, then 
$G\setminus S$ is sequentially Cohen-Macaulay. 
\end{corollary}

\begin{proof} We may assume that $S=\{x_1,\ldots,x_s\}$. Set 
$G_0=G$ and $G_i=G_{i-1}\setminus(y_i\cup N_G(y_i))$ for
$i=1,\ldots,s$ where $y_i$ is
the degree 1 vertex adjacent to $x_i$.  
 Notice that $G_s=G\setminus S$. Hence by 
Theorem~\ref{gscm->cscm1} the graph $G\setminus S$ is 
sequentially Cohen-Macaulay. \end{proof}

We make use of the following result of 
Herzog and Hibi that links the notions of componentwise
linearity and sequentially Cohen-Macaulayness.  We begin
by recalling:

\begin{definition}
\label{d:CWL}
Let $(I_d)$ 
denote the ideal generated by all degree $d$ elements of
a homogeneous ideal $I$.
Then $I$ is called {\it componentwise linear} if $(I_d)$ has a linear
resolution for all $d$.
\end{definition}

\begin{definition}  If $I$ is a squarefree monomial ideal, then
the {\it squarefree Alexander dual} of
$I = (x_{1,1}\cdots x_{1,{s_1}},\ldots,x_{t,1}\cdots x_{t,{s_t}})$
is the ideal
\[I^{\vee} =
(x_{1,1},\ldots,x_{1,s_1}) \cap \cdots \cap (x_{t,1},\ldots,x_{t,s_t}).\]
\end{definition}

If $I$ is a square-free
monomial ideal we write $I_{[d]}$ for the ideal generated by all the
squarefree monomial ideals of degree $d$ in $I$.  
 
\begin{theorem}{\rm (\cite{HH})} \label{t.seqcm}
Let $I$ be a squarefree monomial ideal of $R$. Then 
\begin{itemize}
\item[\rm (a)] $R/I$ is sequentially Cohen-Macaulay if and only if $I^{\vee}$
is componentwise linear.

\item[\rm (b)] $I$ is componentwise linear if and only if 
$I_{[d]}$ has a linear resolution for all $d \geq 0$.
\end{itemize}
\end{theorem}

\begin{lemma}\label{oct15-06} Let $G$ be a bipartite graph with 
bipartition $\{x_1,\ldots,x_m\}$, $\{y_1,\ldots,y_n\}$. If $G$ is
sequentially Cohen-Macaulay, then there is $v\in V_G$ with $\deg(v)=1$.
\end{lemma}

\begin{proof} We may assume that $m\leq n$ and that $G$ has no isolated
vertices. Let $J$ be the Alexander dual of $I=I(G)$ and 
let $L=J_{[n]}$ be the monomial ideal generated by the 
square-free monomials of $J$ of degree $n$. We may assume that $L$ is
generated by $g_1,\ldots,g_q$, where $g_1=y_1y_2\cdots y_n$ and 
$g_2=x_1\cdots x_my_1\cdots y_{n-m}$. Consider the linear map
$$
R^q\stackrel{\varphi}{\longrightarrow} R\ \ \ (e_i\mapsto g_i).
$$
The kernel of this map is generated by syzygies of the 
form
$$
(g_j/\gcd(g_i,g_j))e_i-(g_i/\gcd(g_i,g_j))e_j.
$$   
Since the vector $\alpha=x_1\cdots x_me_1-y_{n-m+1}\cdots y_ne_2$ is
in ${\rm ker}(\varphi)$ and since ${\rm ker}(\varphi)$ is generated by
linear syzygies (see Theorem~\ref{t.seqcm}), there is a linear syzygy
of $L$  of the form  
$x_je_1-ze_k$, where $z$ is a variable, $k\neq 1$. Hence 
$x_j(y_1\cdots y_n)=z(g_k)$ and $g_k=x_jy_1\cdots
y_{i-1}y_{i+1}\cdots y_n$ for some $i$.  Because the support of $g_k$ is a
vertex cover of $G$, we get that the complement of the support of 
$g_k$, i.e., $\{y_i,x_1,\ldots,x_{j-1},x_{j+1},\ldots,x_m\}$, is an
independent set of $G$. Thus $y_i$ can only be adjacent to $x_j$, 
i.e., $\deg(y_i)=1$. \end{proof}

We come to the main result of this section.

\begin{theorem} \label{SCM=shellable}
Let $G$ be a bipartite graph. Then $G$ is
shellable if and only if $G$ is sequentially Cohen-Macaulay. 
\end{theorem}

\begin{proof} Since $G$ shellable implies $G$ 
sequentially Cohen-Macaulay (Theorem \ref{shellable->scm}) we only need
to show the converse. Assume that $G$ is sequentially
Cohen-Macaulay. The proof is by induction on the number of
vertices of $G$. By Lemma~\ref{oct15-06} there is a vertex $x_1$ of
$G$ of degree $1$. Let $y_1$ be the vertex of $G$ adjacent to $x_1$.
Consider the subgraphs $G_1=G\setminus(\{x_1\}\cup N_G(x_1))$ and 
$G_2=G\setminus(\{y_1\}\cup N_G(y_1))$. By 
Theorem~\ref{gscm->cscm1} $G_1$ and $G_2$ are sequentially
Cohen-Macaulay. Hence $\Delta_{G_1}$ and $\Delta_{G_2}$ are shellable
by the induction hypothesis. Therefore $\Delta_G$ is shellable by 
Theorem~\ref{oct15-1-06}. \end{proof}

As we saw in Corollary \ref{recursivebuild}, 
one can verify recursively that a bipartite graph
is shellable.  The above theorem, therefore, implies the same
for sequentially Cohen-Macaulay bipartite graphs.
In particular, we have:

\begin{corollary}\label{scm-build}
Let $G$ be a bipartite graph.  Then $G$ is sequentially Cohen-Macaulay if and only
if there are adjacent vertices $x$ and $y$ 
with $\deg(x)=1$ such that the bipartite graphs $G \setminus (\{x\} \cup N_G(x))$
and $G \setminus (\{y\} \cup N_G(y))$ are sequentially Cohen-Macaulay.
\end{corollary}

\begin{example}
No even cycle $\mathcal{C}_{2m}$ can be sequentially Cohen-Macaulay
since $\mathcal{C}_{2m}$ is 
a bipartite graph that
does not have a vertex of degree 1.
\end{example}

%%%%%%%%%%%%%%%%%%%%%%%%%%%%%%%%%%%%%%%%%%%%%%%%%%%%%%%%%%%%%%%%%%%

\section{An application to Cohen-Macaulay bipartite graphs}\label{s.main}

If $G$ is a bipartite graph without isolated vertices whose edge
ideal $I(G)$ is unmixed, 
i.e., all the associated primes of $I(G)$ have the same height, then one can
show (see, for example, \cite[Theorem 6.4.2]{V}) 
that $G$ must have the following two properties: 
\begin{enumerate}
\item[$(1)$]
if $V_1 = \{x_1,\ldots,x_g\}$ and $V_2 = \{y_1,\ldots,y_h\}$ 
is the bipartition of $V_G$, then $g = h$, 
\item[$(2)$] for $i=1,\ldots,g$, (after relabeling) $\{x_i,y_i\}$ is an edge
of $G$.  
\end{enumerate}
Properties $(1)$ and $(2)$ are deduced from
the fact that all the minimal vertex covers of a graph
whose edge ideal is unmixed ideal must have the same size.  Cohen-Macaulay
bipartite graphs are, therefore, a subset of all the graphs that 
satisfies $(1)$ and $(2)$ since their edge ideals are unmixed.

If $G$ is any bipartite graph that satisfies $(1)$ and $(2)$,
then Carr\'a Ferro and Ferrarello \cite{carra-ferrarello}
introduced a way to construct a directed graph from
the graph $G$.  Precisely,  we define a directed graph $\mathcal D$ with
vertex set $V_1$ as follows: $(x_i,x_j)$ is a directed edge of $\mathcal D$ 
if $i\neq j$ and $\{x_i,y_j\}$ is an edge of $G$. 
In this section $G$ will be any bipartite graph that satisfies
conditions $(1)$ and $(2)$.  We will show how $G$
being (sequentially) Cohen-Macaulay affects the graph $\mathcal D$.
In particular, we can express Herzog and Hibi's \cite{herzog-hibi}
classification of Cohen-Macaulay bipartite graphs
in terms of the graph $\mathcal D$.

We say that a cycle
$\mathcal{C}$ of ${\mathcal D}$ is {\it oriented\/} if all the arrows
of $\mathcal{C}$ are oriented in the same direction. 

\begin{example}
If $G  = \mathcal{C}_4$ with edge set $\{\{x_1,y_1\},\{x_2,y_2\},\{x_1,y_2\},
\{x_2,y_1\}\}$,
then $\mathcal D$ has two vertices $x_1,x_2$ and
two arrows $(x_1,x_2)$, $(x_2,x_1)$ forming an oriented 
cycle of length 
two. 
\end{example}

\begin{lemma}{\rm \cite[Theorem~16.3(4),
p.~200]{Har}}\label{acyclic-char}  Let $\mathcal D$ be the directed
graph described above.
$\mathcal D$ is {\it acyclic}, 
i.e., $\mathcal D$ has no oriented cycles,
if and only if there is a linear  
ordering of the vertex set $V_1$ such that all the 
edges of $\mathcal D$ are of the form 
$(x_i,x_j)$ with $i<j$. 
\end{lemma}

Recall that $\mathcal D$ is called {\it
transitive} if for any two $(x_i,x_j)$,
$(x_j,x_k)$ in $E_{\mathcal D}$ with $i,j,k$ distinct, we have 
that $(x_i,x_k)\in E_{\mathcal D}$. 

\begin{theorem}{\rm(\cite{unmixed})}\label{unmixed-char} 
Let $G$ be a bipartite graph satisfying $(1)$ and $(2)$.
The digraph $\mathcal D$ is 
transitive if and only if $G$ is unmixed, i.e., all minimal vertex covers
of $G$ have the same cardinality. 
\end{theorem}

We can now show how $G$ being sequentially Cohen-Macaulay affects the graph $\D$.

\begin{theorem}\label{chordalCWL}
Let $G$ be a bipartite graph satisfying $(1)$ and $(2)$.
If $G$ is sequentially Cohen-Macaulay, then the directed graph $\D$ is acyclic. 
\end{theorem}

\begin{proof} We proceed by induction on the number of 
vertices of $G$. Assume that $\mathcal D$ has an oriented cycle 
$\mathcal{C}_r$ with vertices $\{x_{i_1},\ldots,x_{i_r}\}$. This means that the
graph $G$ has a cycle 
$$
\mathcal{C}_{2r}=\{y_{i_1},x_{i_1},y_{i_2},x_{i_2},y_{i_3},
\ldots,y_{i_{r-1}},x_{i_{r-1}},
y_{i_r},x_{i_r}\}
$$
of length $2r$. By Lemma~\ref{oct15-06}, the graph $G$ has a vertex
$v$ of degree $1$. Notice that 
$v\notin\{x_{i_1},\ldots,x_{i_r},y_{i_1},\ldots,y_{i_r}\}$. 
Furthermore, if $w$ is the vertex adjacent to $v$, we also have
$w \notin \{x_{i_1},\ldots,x_{i_r},y_{i_1},\ldots,y_{i_r}\}$.  Hence 
by Theorem~\ref{gscm->cscm1} the graph $G'=G\setminus(\{v\}\cup
N_G(v))$ is sequentially Cohen-Macaulay and ${\mathcal D}_{G'}$ has
an oriented cycle, a 
contradiction to the induction hypotheses.
 Thus $\mathcal D$ has no oriented cycles, as required. 
\end{proof}

\begin{example} The converse of the above theorem does not hold
as illustrated through the following example.
Let $G$ be the graph

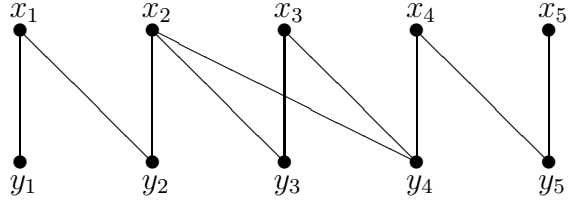
\begin{picture}(100,60)(-100,30)
\put(0,0){\line(0,1){50}}
\put(0,0){\circle*{5}}
\put(0,50){\circle*{5}}
\put(0,50){\line(1,-1){50}}
\put(50,50){\line(1,-1){50}}
\put(50,50){\line(2,-1){100}}
\put(100,50){\line(1,-1){50}}
\put(50,0){\line(0,1){50}}
\put(100,0){\line(0,1){50}}
\put(150,0){\line(0,1){50}}
\put(200,0){\line(0,1){50}}
\put(150,50){\line(1,-1){50}}
 
\put(50,0){\circle*{5}}
\put(100,0){\circle*{5}}
\put(150,0){\circle*{5}}
\put(200,0){\circle*{5}}
\put(50,50){\circle*{5}}
\put(100,50){\circle*{5}}
\put(150,50){\circle*{5}}
\put(200,50){\circle*{5}}

\put(-4,55){$x_1$}
\put(46,55){$x_2$}
\put(96,55){$x_3$}
\put(146,55){$x_4$}
\put(196,55){$x_5$}

\put(-4,-10){$y_1$}
\put(46,-10){$y_2$}
\put(96,-10){$y_3$}
\put(146,-10){$y_4$}
\put(196,-10){$y_5$}
\end{picture}
\vspace{2cm}

\noindent
By Lemma \ref{acyclic-char} $G$ is a bipartite graph whose directed graph
$\mathcal{D}$ is acyclic.  However, $G$ is not sequentially
Cohen-Macaulay.  To verify this, 
note that by Corollary \ref{scm-build}, if $G$ is sequentially Cohen-Macaulay,
then $G_1 = G\setminus(\{x_5\} \cup N_G(x_5))$ and $G_2 =
G\setminus(\{y_5\} \cup N_G(y_5))$ 
are sequentially Cohen-Macaulay.  (Note that by the symmetry of the
graph, we can use 
either $\{x_5,y_5\}$ or $\{x_1,y_1\}$.)
But $G_2$ is sequentially Cohen-Macaulay if and only if $H_1 = G_2
\setminus(\{y_1\}\cup N_G(y_1))$ 
and $H_2 = G_2 \setminus(\{x_1\} \cup N_G(x_1))$ are sequentially Cohen-Macaulay.
But $H_2$ is the graph of $\mathcal{C}_4$ which is not sequentially Cohen-Macaulay.  
Hence, $G$ is not sequentially Cohen-Macaulay.
\end{example}

Bipartite Cohen-Macaulay graphs have been studied in
\cite{EV,herzog-hibi,V}. In \cite{EV} it is shown that $G$ is a
Cohen-Macaulay bipartite graph if and only if $\Delta_G$ is pure 
shellable. In \cite{herzog-hibi} Herzog and Hibi give a graph
theoretical description of  
Cohen-Macaulay bipartite graphs. This description 
can be expressed in terms of $\mathcal D$, as was pointed out in 
\cite{carra-ferrarello}. As a corollary, we prove Herzog and Hibi's
result classifying Cohen-Macaulay bipartite graphs.

\begin{corollary}{\rm(\cite{carra-ferrarello,herzog-hibi})}\label{c.hhz}
Let $G$ be a bipartite graph satisfying $(1)$ and $(2)$.
Then
$G$ is Cohen-Macaulay if and only
if $\mathcal D$ is acyclic and transitive. 
\end{corollary} 

\begin{proof}
($\Rightarrow$) By Theorem~\ref{unmixed-char}, $\mathcal D$ is
transitive, and by Theorem~\ref{chordalCWL}, $\mathcal D$ is
acyclic. 

($\Leftarrow$) The proof is by induction on $g = |V_1|$. The case $g=1$ is
clear. We may assume that $G$ is connected and $g\geq 2$. By
Lemma~\ref{acyclic-char} we may also assume that 
if $\{x_i,y_j\} \in E_G$, then $i \leq j$. Let 
$N_G(y_g)=\{x_{r_1},\ldots,x_{r_s}\}$ be the set of all 
vertices of $G$ adjacent to $y_g$, where $x_{r_s}=x_g$. 
Consider the subgraph
$G'=G\setminus (\{y_g\}\cup N_G(y_g))$. We claim that
$y_{r_1},\ldots,y_{r_{s-1}}$ are isolated vertices of $G'$. Indeed if
$y_{r_j}$ is not isolated, there is an edge $\{x_i,y_{r_j}\}$ in $G'$
with $i<r_j$. Hence, by the transitivity of $\mathcal D$, we get that
$\{x_i,y_g\}$ is an edge 
of $G$ and $x_i$ must be a vertex in $N_G(y_g)$, a contradiction. Thus, by
induction, the graphs $G'$ and 
$G''= G\setminus\{x_g,y_g\}=G\setminus(\{x_g\}\cup N_G(x_g))$ are
Cohen-Macaulay.   If $R_1 = k[x ~|~ x \in V_{G'}]$ and
$R_2 = k[x ~|~ x \in V_{G''}]$, then by induction
$\dim R_1/I(G') = g-s$ and $\dim R_2/I(G'') = g-1$.
Since $(I(G)\colon y_g) = (x_{i_1},\ldots,x_{r_s},I(G'))$ and
$(y_g,I(G)) = (y_g,I(G''))$,  the 
ends of the exact sequence
$$
0\longrightarrow R/(I(G)\colon y_g)\stackrel{y_g}{\longrightarrow} 
R/I(G)\longrightarrow R/(I(G),y_g)\longrightarrow 0
$$
are Cohen-Macaulay modules of dimension $g$.
On the other hand, because $\mathcal D$ is transitive,
by Theorem \ref{acyclic-char} the graph $G$ is unmixed,
and thus $\dim R/I(G) = 
\dim R - {\rm ht}(I(G)) = 2g - g = g$ since $g$ 
is the size of any minimal vertex covering.  Consequently by applying
the depth lemma (see \cite[Corollary 18.6]{E}) to the above
short exact sequence, we have 
\[\dim R/I(G) \geq \operatorname{depth} R/I(G) \geq
\min\{\operatorname{depth}~R/(I(G)\colon y_g),\operatorname{depth}~R/(I(G), y_g)+1\} = g\]
whence $R/I(G)$  is Cohen-Macaulay of dimension $g$. 
\end{proof}

Cohen-Macaulay trees, first studied in \cite{Vi2}, can also be described in terms of $\mathcal D$:

\begin{theorem}  Let $G$ be a tree satisfying $(1)$ and $(2)$.
Then $G$ is a Cohen-Macaulay tree if and only if $\mathcal D$ is a 
tree such that every vertex $x_i$ of $\mathcal D$ is either a {\it
source\/} {\rm (}i.e., has only arrows leaving $x_i${\rm )} or a
{\it sink\/} 
{\rm (}i.e., 
has only arrows entering $x_i${\rm )}. 
\end{theorem}

\begin{proof}
($\Rightarrow$) 
Since a tree is bipartite, $\mathcal{D}$ is both acyclic and transitive.
Suppose there is a vertex $x_i$ that is not a sink or source. i.e.,
there is an arrow entering $x_i$ and one leaving $x_i$.  Suppose
the arrow entering $x_i$ originates at $x_j$, and the arrow leaving
$x_i$ goes to $x_k$.  Note that $x_j \neq x_k$ because
otherwise we would have a cycle in the acyclic graph $\mathcal{D}$.
Because $\mathcal{D}$ is transitive, the directed edge $(x_j,x_k)$
also belongs to $\mathcal{D}$.  But then the induced graph
on the vertices $\{x_j,y_i,x_i,y_k\}$ in $G$ forms the cycle $\mathcal{C}_4$,
contradicting the fact that $G$ is a tree.

$(\Leftarrow)$ The hypotheses on $\mathcal{D}$ imply $\mathcal{D}$
is acyclic and transitive, so apply Theorem \ref{c.hhz}.
\end{proof}
%%%%%%%%%%%%%%%%%%%%%%%%%%%%%%%%%%%%%%%%%%%%%%%%%%%%%%%%%%%%%%%%%%%

\section{Clutters with the free vertex property are shellable}

We now extend the scope of our paper to include a special
family of hypergraphs called clutters.  The results of this section
allow us to give a new proof to a result of Faridi on the
sequentially Cohen-Macaulayness 
of simplicial forests.
 
A {\it clutter\/} $\mathcal C$ with 
vertex set $X=\{x_1,\ldots,x_n\}$ is a family of subsets of $X$,
called edges, none of which is included in another. The set of
vertices and edges of $\mathcal C$ are denoted by $V_{\mathcal C}$ and
$E_{\mathcal C}$ respectively. A basic example  
of a clutter is a graph.  Note that a clutter is an example of a hypergraph
on the vertex set of $X$;  a clutter is sometimes called a simple hypergraph,
as in \cite{HVT}. For a thorough study of clutters---that includes 18 
conjectures in the area---from the point of
view of combinatorial optimization see \cite{cornu-book}. 

Let $R=k[x_1,\ldots,x_n]$ be a polynomial ring 
over a field $k$ and let $I$ be an ideal 
of $R$ minimally generated by a finite set
$\{x^{v_1},\ldots,x^{v_q}\}$ of
square-free monomials. As usual we
use  $x^a$ as an abbreviation for $x_1^{a_1} \cdots x_n^{a_n}$, 
where $a=(a_1,\ldots,a_n)\in \mathbb{N}^n$. Note that the entries of
each $v_i$ are in $\{0,1\}$. We associate to the 
ideal $I$ a {\it clutter\/} $\mathcal C$ by taking the set 
of indeterminates $V_{\mathcal C}=\{x_1,\ldots,x_n\}$ as the vertex set and 
$E_{\mathcal C}=\{S_1,\ldots,S_q\}$ as the edge set, where 
$S_i={\rm supp}(x^{v_i})$ is the {\it support\/} of 
$x^{v_i}$, i.e., $S_i$ is the set of variables that occur in
$x^{v_i}$.  For this
reason  $I$ is called the {\it edge ideal\/} of $\mathcal C$ and is denoted
$I = I(\mathcal C)$.   Edge ideals of clutters
are also called {\it facet ideals} \cite{Faridi} because 
$S_1,\ldots,S_q$ are exactly the facets of the simplicial complex
$\Delta=\langle S_1,\ldots,S_q\rangle$ generated by $S_1,\ldots,S_q$.

A subset $C\subset X$ is a 
{\it minimal vertex cover\/} of the clutter $\mathcal C$ if: 
(i) every edge of $\mathcal C$ contains at least one vertex of $C$, 
and (ii) there is no proper subset of $C$ with the first 
property. If $C$ only satisfies condition (i), then $C$ is 
called a {\it vertex cover\/} of $\mathcal C$. Notice that
$\mathfrak{p}$ is a minimal prime of $I =I(\mathcal C)$ if and only if 
$\mathfrak{p}=(C)$ for some minimal vertex cover $C$ of $\mathcal C$.
In particular, if $D_1,\ldots,D_t$ is a complete list
of the minimal vertex covers of $\mathcal C$, then
\[I(\mathcal C) =(D_1)\cap (D_2)\cap\cdots\cap(D_t).\]

Because $I = I(\mathcal C)$ is a square-free monomial ideal,
it also corresponds to a simplicial complex via 
the Stanley-Reisner correspondence \cite{Stanley}.  
We let $\Delta_{\mathcal C}$ represent 
this simplicial complex. 
Note that $F$ is a facet of $\Delta_{\mathcal C}$ if and only if 
$X\setminus F$ is a minimal vertex cover of $\mathcal C$. 
As for graphs, we may say 
that the clutter $\mathcal C$ is {\it shellable\/} if
$\Delta_{\mathcal C}$ is shellable.

\begin{lemma}\label{induced-shelling} Let $\mathcal{C}$ be a clutter
with minimal vertex covers $D_1,\ldots,D_t$.
If $\Delta_{\mathcal{C}}$ is shellable and $A\subset
V_{\mathcal{C}}$ is a set
of vertices, then the Stanley-Reisner complex $\Delta_{I'}$ of 
the ideal 
$$
I'=\bigcap_{\scriptstyle D_i\cap A=\emptyset}(D_i)
$$
is shellable with respect to the linear ordering of the facets of
$\Delta_{I'}$ induced by the shelling of the simplicial complex
$\Delta_{\mathcal{C}}$. 
\end{lemma}

\begin{proof} Let $H_1,\ldots,H_t$ be a shelling of $\Delta_{\mathcal{C}}$. We may
assume that $H_i=V_{\mathcal{C}}\setminus D_i$ for all $i$. Let $H_i$ and $H_j$ be
two facets of $\Delta_{I'}$ with
$i<j$, i.e., $A\cap D_i=\emptyset$ and $A\cap D_j=\emptyset$. 
By the shellability of $\Delta_{\mathcal{C}}$, there 
is an $x\in H_j\setminus H_i$ and an $\ell<j$ such that $H_j\setminus
H_\ell=\{x\}$. It suffices to prove that
$D_\ell\cap A=\emptyset$. If $D_\ell\cap A\neq\emptyset$, pick $z\in
D_\ell \cap A$. Then $z\notin D_i\cup D_j$ and $z\in H_i\cap H_j$.
Since $z\notin H_\ell$ (otherwise  $z\notin D_\ell$, a contradiction),
we get $z\in H_j\setminus H_\ell$, i.e., $z=x$, a contradiction
because $x\notin H_i$. 
\end{proof}

An ideal $I'$ is called a {\it minor\/} of $I$ if there is a subset 
$X'=\{x_{i_1},\ldots,x_{i_r},x_{j_1},\ldots,x_{j_s}\}$ of the set of
variables $X=\{x_1,\ldots,x_n\}$ such that $I'$ is a proper ideal of 
$R'=k[X\setminus X']$ that can be obtained from a generating set of
$I$ by setting 
$x_{i_k}=0$ and $x_{j_\ell}=1$ for all $k,\ell$. The ideal $I$ is
also considered to be a minor.  A {\it minor\/} of $\mathcal C$ 
is a clutter ${\mathcal C}'$ on the vertex set 
$V_{\mathcal{C}'}=X\setminus X'$ 
that corresponds to a minor $(0)\subsetneq
I'\subsetneq R'$. Notice that the edges of ${\mathcal C}'$ are obtained from $I'$
by considering the unique set  
of square-free monomials of $R'$ that minimally generate $I'$.   
For use below we say $x_i$ is a {\it free variable\/} (resp. {\it free
vertex}) of $I$ (resp. ${\mathcal C}$) if $x_i$
only appears in one of the monomials $x^{v_1},\ldots,x^{v_q}$ (resp. in
one of the edges of $\mathcal C$). If all the minors of 
$\mathcal C$ have free vertices, we say that ${\mathcal C}$ 
has the {\it free vertex property}.   Note that if $\mathcal C$
has the free vertex property, then so do all of its minors.

\begin{lemma}\label{covers-xn} Let $x_n$ be a free variable of $I = I(\mathcal{C})
= (x^{v_1},\ldots,x^{v_{q-1}},x^{v_q})$,
and let $x^{v_q}=x_nx^u$. {\rm (a)} If $\mathcal{C}_1$ is the clutter associated
to $J=(x^{v_1},\ldots,x^{v_{q-1}})$, then $C$ is a minimal vertex
cover of $\mathcal{C}$ containing $x_n$ if and only if 
$C\cap{\rm supp}(x^{u})=\emptyset$ and $C=\{x_n\}\cup C'$ for some
minimal vertex cover $C'$ of $\mathcal{C}_1$. {\rm (b)} If $\mathcal{C}_2$
is the clutter associated 
to $L=(x^{v_1},\ldots,x^{v_{q-1}},x^u)$, then $C$ is a minimal vertex
cover of $\mathcal{C}$ not containing $x_n$ if and only if $C$ is a
minimal vertex cover of $\mathcal{C}_2$.
\end{lemma}

\begin{proof} (a) Assume that $C$ is a minimal vertex cover of
$\mathcal{C}$ containing $x_n$. If $C\cap{\rm
supp}(x^{u})\neq\emptyset$, then $C\setminus\{x_n\}$ is a vertex cover
of $\mathcal{C}$, a contradiction. Thus $C\cap{\rm
supp}(x^{u})=\emptyset$. Hence it suffices to notice that
$C'=C\setminus\{x_n\}$ is a minimal
vertex cover of $\mathcal{C}_1$. The converse also follows readily.

(b) Assume that $C$ is a minimal vertex cover of
$\mathcal{C}$ not containing $x_n$. Let $x^a$ be a minimal generator
of $I(\mathcal{C}_2)$, then either $x^u$ divides $x^a$ or
$x^a=x^{v_i}$ for some $i<q$. Then clearly $C\cap {\rm
supp}(x^a)\neq\emptyset$ because $C\cap A\neq\emptyset$, where 
$A={\rm supp}(x^u)$. Thus $C$ is a vertex cover of $\mathcal{C}_2$. 
To prove that $C$ is minimal take $C'\subsetneq C$. We must show that
there is an edge of $\mathcal{C}_2$ not covered by $C'$. As $C$ is a
minimal vertex cover of $\mathcal{C}$, there is $x^{v_i}$ such 
that ${\rm supp}(x^{v_i})\cap C'=\emptyset$. If $x^{v_i}$ is a minimal
generator of $\mathcal{C}_2$ there is nothing to prove, otherwise
$x^u$ divides $x^{v_i}$ and the edge $A$ of $\mathcal{C}_2$ is not 
covered by $C'$. The converse also follows readily.
\end{proof}

\begin{theorem} \label{fvp}
If the clutter $\mathcal{C}$ has the free vertex property, then
$\Delta_{\mathcal{C}}$ is shellable.
\end{theorem}

\begin{proof} We proceed by induction on the number of vertices of 
$\mathcal C$. Let $x_n$ be a free variable of $I=I(\mathcal{C}) =
(x^{v_1},\ldots,x^{v_{q-1}},x^{v_q})$. We
may assume that $x_n$ occurs in $x^{v_q}$. Hence we can write
$x^{v_q}=x_nx^u$ for some $x^u$ such that 
$x_n\notin{\rm supp}(x^u)$. For use below we set $A={\rm supp}(x^u)$.
Consider the ideals $J=(x^{v_1},\ldots,x^{v_{q-1}})$ and
$L=(J,x^u)$. Then 
$J=I({\mathcal C}_1)$ and $L=I({\mathcal C}_2)$, 
where ${\mathcal C}_1$ and ${\mathcal C}_2$ are the clutters 
defined by the ideals $J$ and $L$, respectively. 
Notice that $J$ and $L$ are minors of 
the ideal $I$ obtained by setting $x_n=0$  
and $x_n=1$, respectively. The vertex set of $\mathcal{C}_i$ is 
$V_{\mathcal{C}_i}=X\setminus\{x_n\}$ for $i=1,2$. Thus $\Delta_{\mathcal{C}_1}$ and
$\Delta_{\mathcal{C}_2}$ are  
shellable by the induction hypothesis. Let $F_1,\ldots,F_r$ be the facets
of $\Delta_{\mathcal C}$ that contain $x_n$ and let $G_1,\ldots,G_s$ be the facets
of $\Delta_{\mathcal C}$ that do not contain $x_n$. Set
$C_i=X\setminus G_i$ and $C_i'=C_i\setminus \{x_n\}$ for
$i=1,\ldots,s$. Then $C_1,\ldots,C_s$ is the set of minimal vertex
covers of $\mathcal{C}$ that contain $x_n$, and by
Lemma~\ref{covers-xn}(a) $C_1',\ldots,C_s'$ is the set of minimal vertex
covers of $\mathcal{C}_1$ that do not intersect $A$. One has the
equality $G_i=V_{\mathcal{C}_1}\setminus C_i'$ for all $i$. Hence, by 
the shellability of $\Delta_{\mathcal{C}_1}$ and using
Lemma~\ref{induced-shelling}, we may assume that 
$G_1,\ldots,G_s$ is a shelling for the simplicial complex generated by
$G_1,\ldots,G_s$. By Lemma~\ref{covers-xn}(b) one has that $C$ is a minimal vertex
cover of $\mathcal{C}$ not containing $x_n$ if and only if $C$ is a
minimal vertex cover of $\mathcal{C}_2$.  Thus, $F$ is a facet of $\Delta_{\mathcal{C}}$
that contains $x_n$, i.e., $F = F' \cup \{x_n\}$ if and only if 
$F'$ is a facet of  $\Delta_{\mathcal{C}_2}$.
By induction we may also assume that
$F'_1=F_1\setminus\{x_n\},\ldots,F'_r=F_r\setminus\{x_n\}$ 
is a shelling of $\Delta_{\mathcal{C}_2}$. We now prove that
$$
F_1,\ldots,F_r,G_1,\ldots,G_s ~~\mbox{with $F_i = F'_i \cup \{x_n\}$}
$$
is a shelling of $\Delta_{\mathcal{C}}$. We need only show that 
given $G_j$ and $F_i$ there is $a\in G_j\setminus F_i$ and $F_\ell$
such that $G_j\setminus F_\ell=\{a\}$. We can write 
$$
G_j=X\setminus C_j\ \mbox{ and }\ F_i=X\setminus C_i,
$$
where $C_j$ (resp. $C_i$) is a minimal vertex cover of $\mathcal C$
containing $x_n$ (resp. not containing $x_n$). Recall that 
$A={\rm supp}(x^u)$ is an edge of $\mathcal{C}_2$. Notice the
following: (i) $C_j=C_j'\cup\{x_n\}$ for some minimal vertex cover $C_j'$ 
of ${\mathcal C}_1$ such that $A\cap C_j'=\emptyset$, and (ii) $C_i$
is a minimal vertex cover of 
${\mathcal C}_2$. From (i) we get that $A\subset G_j$. Observe that
$A\not\subset F_i$, otherwise $A\cap C_i=\emptyset$, a contradiction
because $C_i$ must cover the edge $A={\rm supp}(u)$. Hence there is
$a\in A\setminus F_i$ and $a\in G_j\setminus F_i$.  Since 
$C_j'\cup\{a\}$ is a vertex cover of $\mathcal{C}$, there is 
a minimal vertex cover $C_\ell$ of $\mathcal{C}$ contained in
$C_j'\cup\{a\}$. Clearly $a\in C_\ell$ because $C_\ell$ has to cover
$x^u$ and $C_j'\cap A=\emptyset$.  Thus
$F_\ell=X\setminus C_\ell$ is a facet of $\Delta_{\mathcal{C}}$ containing $x_n$.
To finish the proof we now prove that $G_j\setminus F_\ell=\{a\}$. We
know that $a\in G_j$. If $a\in F_\ell$, then $a\notin C_\ell$, a
contradiction. Thus $a\in G_j\setminus F_\ell$. Conversely take $z\in
G_j\setminus F_\ell$. Then $z\notin C_j'\cup\{x_n\}$ and 
$z\in C_\ell\subset C_j'\cup\{a\}$. Hence $z=a$, as required. \end{proof}

The $n\times q$ matrix $A$ with column vectors
$v_1,\ldots,v_q$ is called the {\it incidence
matrix\/} of $\mathcal C$. This matrix has entries in $\{0,1\}$. We
say that $A$ (resp. $\mathcal{C}$) is a {\it totally balanced\/}
matrix (resp. clutter) if $A$ has no square 
submatrix of order at least $3$ with exactly two $1$'s in 
each row and column. According to \cite[Corollary~83.3a]{Schr2} a
totally balanced clutter satisfies the free vertex property. Thus we
obtain:

\begin{corollary} If $\mathcal{C}$ is a totally balanced clutter,
then $\Delta_{\mathcal{C}}$ is shellable.
\end{corollary} 
 
Faridi \cite{Faridi} introduced the notion of a leaf for a simplicial
complex $\Delta$.  Precisely, a facet $F$ of $\Delta$ is a {\it leaf} if $F$ is
the only facet of $\Delta$, or there exists a facet $G \neq F$ in $\Delta$ such that
$F \cap F' \subset F \cap G$ for all facets $F' \neq F$ in $\Delta$.
A simplicial complex $\Delta$ is a {\it simplicial forest} if every nonempty
subcollection, i.e., a subcomplex whose
facets are also facets of $\Delta$, of $\Delta$ contains a leaf. 
We can translate Faridi's definition into hypergraph language; we call the
translated version of Faridi's leaf a $f$-leaf.

\begin{definition} An edge $E$ of a clutter $\mathcal{C}$ is aa {\it $f$-leaf} if $E$
is the only edge of $\mathcal{C}$, or if there exists an edge $H$ of
$\mathcal{C}$ such that 
$E \cap E' \subset E \cap H$ for all edges $E' \neq E$ of $\mathcal{C}$.
A clutter $\mathcal{C}$ is an $f$-{\it forest}, if every 
subclutter of $\mathcal{C}$, including $\mathcal{C}$ itself, contains 
an $f$-leaf.  
\end{definition}

In \cite[Theorem~3.2]{hhtz} it is shown that $\mathcal{C}$ is an $f$-forest if
and only if $\mathcal{C}$ is a totally balanced clutter. Thus we
obtain: 

\begin{corollary}\label{maintheorem-simplicial} 
If the clutter $\mathcal{C}$ is an $f$-forest, 
then $\Delta_{\mathcal{C}}$ is shellable.  
\end{corollary}

We now recover the main
result of Faridi \cite{Faridi}:

\begin{corollary}  Let $I = I(\Delta)$ be the facet ideal
of a simplicial forest.  Then
$R/I(\Delta)$ is 
sequentially Cohen-Macaulay.
\end{corollary}

\begin{proof}
If $\Delta = \langle F_1,\ldots,F_s\rangle$, then $I(\Delta)$
is also the edge ideal of the clutter $\mathcal{C}$ whose
edge set is $E_{\mathcal{C}} = \{F_1,\ldots,F_s\}$.  Now apply
Corollary \ref{maintheorem-simplicial}
and Theorem \ref{shellable->scm}.
\end{proof}

\begin{remark}  Since submitting this paper, Soleyman Jahan and Zheng
\cite[Theorem 3.4]{SJZ} have given a generalization of Theorem \ref{fvp}
using the notion a {\it pretty clean monomial ideal}.
\end{remark}
\noindent
{\bf Acknowledgments.} We gratefully acknowledge the computer algebra
system {\tt CoCoA} \cite{Co} which was invaluable in our work on this paper.
The first author also acknowledges the financial support of NSERC;
the second author acknowledges the financial support of 
CONACyT grant 49251-F and SNI.  We also thank the referees for their
careful reading of the paper and for the improvements that
they suggested.

\bibliographystyle{plain}

\end{document}